\documentclass[12pt]{article}
\usepackage{amsmath,amsfonts,amssymb}
\usepackage{graphicx,color}

\usepackage{theorem}

\theorembodyfont{\upshape}

\newcommand{\ZZ}{{\mathbb Z}}

\newcommand{\CC}{{\mathbb C}}
\newcommand{\NN}{{\mathbb N}}
\newcommand{\Der}{Der}
\def\gg{{\mathfrak{g}}}

\def\hh{{\mathfrak{h}}}

\numberwithin{equation}{section}

\begin{document}

\begin{center}
{\LARGE{\bf On Integrable modules for the twisted full toroidal Lie algebra }}  \\ [5mm] 
{\bf  S.Eswara Rao, Punita Batra} \\
{School of Mathematics, Tata Institute of Fundamental Research, Mumbai.}\\

{Harish-Chandra Research Institute, Allahabad, 211019, India.} \\
{ senapati@math.tifr.res.in,
batra@hri.res.in
}\\

\end{center}

\begin{abstract}
The paper is to classify irreducible integrable modules for the twisted 
full toroidal Lie algebra with some technical
conditions. The  twisted full toroidal Lie algebra are extensions of
 multiloop algebra twisted by several finite
order automorphisms. This result genaralizes a result by Fu Jiayuan and
 Cuipo Jiang $[FJ]$, where they consider only one automorphism. \\


\end{abstract}

\section*{Introduction} The main purpose of this paper is to classify
 irreducible integrable modules for the twisted full toroidal Lie algebra 
under certain technical conditions. The  twisted full toroidal Lie algebra 
in several variables is defined using several commuting finite order 
automorphisms of the underlying finite dimensional simple Lie algebra 
$\gg$. The result of this paper generalizes the main theorem of $[FJ]$,
 where they consider only one automorphism.\

The twisted full toroidal Lie algebra is a natural genaralization of the 
classical twisted affine Lie algebra. The classical procedure of realizing
twisted affine Lie algebras using loop algebra in one variable proceeds in 
two steps $[K]$. In the first step, the derived algebra modulo its center 
of the affine Lie algebra is constructed as the algebra of a diagram 
automorphism of a finite dimensional simple Lie algebra. In the second step
 the affine Lie algebra is built from the graded loop algebra by forming 
universal central extension (one dimensional center) and adding graded 
algebra derivations (are dimensional).\

The replacement of derived Lie algebra is multiloop algebra (Sec. 1.5) in 
several variables twisted by finitely many automorphisms. Then we consider
the Universal central extension of the multiloop algebra (infinite 
dimensional) and add graded algebra of derivations (infinite dimensional)
which we call twisted full toroidal Lie algebra and denoted by $\tau$.\

In this paper we classify irreducible integrable modules for $\tau$ with 
finite dimensional weight spaces with non-zero 
central action and with some technical conditions that are satisfied for 
the well known Lie algebras called Lie Torus. See [ABFP].\

The contents of the paper are the following. We fix an irreducible 
integrable module for $\tau$ with non-zero central action. 

The central operators act as scalars. Upto choice of co-ordinates (in other
 words upto an automorphism) we can assume that $K_0$ acts as $C_0 > 0$ and 
$K_i$ acts trivially for $i \neq 0$. In the first main Theorem 5.3 we prove 
that $V$ is an highest weight module with some natural triangular 
decomposition 
$$
\tau = \tau^- \oplus \tau^0 \oplus \tau^+.
$$
Let $T$ be the highest weight space  which is  an irreducible 
$\tau_0$-module(Proposition 6.2) and  naturally $\ZZ^n$- graded.  In 
Sections 7, 8 and 9, we describe $T$ as  $\tau_0$-module. We noticed that 
some parts of  $\tau_0$ acts as scalars and hence we only consider a 
subalgebra $L$ (See 7.3) for which $T$ is irreducible.
We  also consider a certain  subalgebra 
$\stackrel{\sim}{L}$  of $L$ (See 8.1)  and consider  a certain finite
 dimensional quotient  $\stackrel{\sim}{V}$ of $T$, which is 
$\stackrel{\sim}{L}$-module. 
Our approach is to describe $\stackrel{\sim}{V}$ as $\stackrel{\sim}{L}$-
module and obtain $T$ as $L$-module.

In the process we define an $L$-module $L(\stackrel{\sim}{V})$ (See 8.4) and 
prove $L(\stackrel{\sim}{V})$ naturally decomposes into $L$-modules ( See 
notes after Lemma 8.9 and 8.10). We prove one of the component is 
isomorphic to $T$. 
We also prove in Theorem 8.17 that $\stackrel{\sim}{V}$ is  irreducible 
$\stackrel{\sim}{L}$-module if and only if $L(\stackrel{\sim}{V})$ 
decomposes into mutually non-isomorphic irreducible $L$-modules. In this 
case $T$ can be identified in a natural way as a component of 
$L(\stackrel{\sim}{V})$. 
In the case $\stackrel{\sim}{V}$ is reducible, we prove that   
 $\stackrel{\sim}{V}$ is  completely reducible  as $\stackrel{\sim}{L}$
-module and all components are isomorphic (See Proposition 8.27).  In this 
case also $T$ is a submodule of $L(\stackrel{\sim}{V})$ but the inclusion 
is twisted.
In the rest of the paper we  describe each component of 
$\stackrel{\sim}{V}$. 
 It turns out to be an irreducible  module for 
$g l_n \oplus \stackrel{\circ}{\gg}$ (Theorem 9.4). See 7.5 for the 
Definition of  $\stackrel{\circ}{\gg}$.  We will indicate in 9.6, given an 
irreducible  module for 
$g l_n \oplus \stackrel{\circ}{\gg}$, how to obtain a module  for  $L$ and 
thereby for $\tau_0$. In the Theorem 9.7 we will state our final result 
that the original module $V$ is an irreducible quotient of an induced 
module of $T$.

\section{Notation and Preliminaries} 

Throughout this paper we will use the following notation

\begin{enumerate}

\item [{(1.1)}] All vector spaces, algebras and tensor products are over 
complex numbers $\CC$. Let $\ZZ, \NN$ and $\ZZ_+$ 
denote integers, non-negative integers and positive integers.

\item [{(1.2)}] Let $\gg$ be a  finite dimensional simple Lie algebra and 
let $(,)$ be a non-degenerate symmetric bilinear form on 
$\gg$. We fix a positive integer $n$. Let $\sigma_0, \sigma_1, \cdots, 
\sigma_n$ be commuting finite order automorphisms of $\gg$ 
of order $m_0, m_1, \cdots, m_n$ respectively. Let $m= ( m_1, \cdots, m_n) 
\in \ZZ^n$. Let $k= ( k_1, \cdots, k_n)$ and $l= ( l_1, \cdots, l_n)$ 
denote vectors in $\ZZ^n$. 

\item [{(1.3)}] Let $\Gamma = m_1 \ZZ \oplus \cdots \oplus m_n \ZZ$ and 
$\Gamma_0 = m_0 \ZZ$. Let $\Lambda= \ZZ^n /\Gamma$ and 

$\Lambda_0= \ZZ /\Gamma_0$. Let $\overline{k}$ and $\overline{l}$ denote the
images in $\Lambda$. For any integers $k_0$ and $l_0$, let $\overline{k}_0$
and $\overline{l}_0$ denote the images in $\Lambda_0$.\\

Let

$$
\begin{array}{lll}
A&= \CC[t_0^{\pm 1}, \cdots, t_n^{\pm 1}],\\
A_n&= \CC[t_1^{\pm 1}, \cdots, t_n^{\pm 1}],\\
A(m)&= \CC[t_1^{\pm m_1}, \cdots, t_n^{\pm m_n}],\\
A(m_0,m)&= \CC[t_0^{\pm m_0}, t_1^{\pm m_1}, \cdots, t_n^{\pm m_n}].\\
\end{array}
$$

\item [{(1.4)}] For $k \in \ZZ^n$, let $t^k= t_1^{k_1} \cdots t_n^{k_n} \in 
A_n$. Let $\Omega A$ be the vector space spanned by 
symbols $t_0^{k_0} t^k K_i, 0 \leq i \leq n, k_0 \in \ZZ, k \in \ZZ^n$. Let 
$dA$ be the subspace spanned by 
$\displaystyle{\sum_{i=0}^n} k_i t_0^{k_0} t^k K_i$.\

Let $L(\gg)= \gg \otimes A$ and notice that it has a natural structure of a 
Lie algebra. We will now define toroidal Lie algebra
$$
\stackrel{\sim}{L}(\gg)= L(\gg) \oplus \Omega A/dA.
$$

Let $X(k_0,k)= X \otimes t_0^{k_0} t^k$ and $Y= Y \otimes t_0^{l_0} t^l$ for 
$X, Y \in \gg, k_0, l_0 \in \ZZ$ and 
$k, l \in \ZZ^n.$

\item [{(1.4.1)}] $[X(k_0,k), Y(l_0,l)= [X,Y] (l_0+k_0, l+k) + (X,Y) 
\sum k_i t_0^{l_0+k_0} t^{l+k} K_i$.

\item [{(1.4.2)}] $\Omega A/dA$ is central.\

It is well known that $\stackrel{\sim}{L} (\gg)$ is the universal central 
extension  of $L(\gg)$. See $[EMY]$ and $[Ka]$.

\item [{(1.5)}] We will now define multiloop algebra as a subalgebra of 
$L(\gg)$. See $[ABFB]$ for more details. For 
$0 \leq i \leq n$, let $\xi_i$ be a $m_i$th primitive root of unity.\\
Let 
$$
\gg(\overline{k}_0, \overline{k})= \{x \in \gg | \sigma_i x = \xi_i^{k_i} x, 
0 \leq i \leq n\}. 
$$

Then define 
$$
L(\gg, {\bf \sigma}) = \displaystyle{\bigoplus_{(k_o,k) \in \ZZ^{n+1}}} 
\gg(\overline{k}_0, \overline{k}) \otimes t_0^{k_0}t^k,
$$

which is called a multiloop algebra.

\item [{(1.6)}] The finite dimensional irreducible modules for $L(\gg, 
\sigma)$ are classified by Michael Lau $[ML]$.

\item [{(1.7)}] Suppose $\hh_1$ is a finite dimensional ad-diagonalizable 
subalgebra of a Lie algebra $\gg_1$. We set for $\alpha \in \hh_1^*$
$$
\gg_{1,\alpha} = \{x \in \gg_1 | [h.x]= \alpha (h) x, h \in \hh_1\}.
$$
Then we have 
$$
\gg_1= \displaystyle{\bigoplus_{\alpha \in \hh^*_1}} \gg_{1,\alpha}.
$$
Let $\Delta(\gg_1, \hh_1) = \{\alpha \in \hh^*_1 | \gg_{1,\alpha} \neq 0\}$,
which includes zero.\\
Let $\Delta^\times (\gg_1, \hh_1) = \Delta (\gg_1, \hh_1)\backslash \{0\}.$

\item [{(1.8)}] We will now define the universal central extension of 
$L(\gg, \sigma).$ Define 
$\Omega A (m_0,m)$ and $dA(m_0,m)$ similar to the definition of $\Omega A$ 
and $dA$ by replacing $A$ by $A(m_0,m)$. Denote $Z(m_0,m) = \Omega A 
(m_0,m)/ dA(m_0,m) $ and note that $Z(m_0,m) \subseteq \Omega A /dA$.\\
Define
$$
\stackrel{\sim}{L} (\gg,\sigma) = L(\gg,\sigma) \oplus Z(m_0,m).
$$
Let $X \in \gg (\overline{k}_0, \overline{k})$ and $Y \in \gg 
(\overline{l}_0, \overline{l})$ and let 
$X(k_0,k) = X \otimes t_0^{k_0} t^k$ and $Y(l_0,l)= Y \otimes t_0^{l_0} t^l$. 
Define 

\item [{(1.8.1)}] $[X(k_0,k), Y(l_0,l)]= [X,Y](k_0+l_0, k+l) + (X,Y) 
\sum k_i t_0^{l_0+k_0} t^{l+k} K_i$.

\item [{(1.8.2)}] $Z(m_0,m)$ is central. \

Notice that $(X,Y) \neq 0 \Rightarrow k+l \in \Gamma$ 
and $k_0+l_0 \in \Gamma_0$. This follows from the standard fact that $(,)$ 
is invariant under $\sigma_i$ for 
$0 \leq i \leq n$. This proves that the above Lie bracket is well defined. 
This Lie bracket is nothing but the 
restriction defined in (1.4).

\item [{(1.9)}] {\bf Proposition}  $\stackrel{\sim}{L} (\gg,\sigma)$ is the 
universal central extension of $L(\gg,\sigma)$. See Corollary (3.27) of [JS].

\end{enumerate}

\section{Derivation algebra of $A(m_0,m)$ and its extension 
to $Z(m_0,m)$. }

\begin{enumerate}

\item [{(2.1)}] Let $D(m_0,m)$ be the derivation algebra of $A(m_0,m)$. 
From now onwards we let $s$ and $r$ to be 
in $\Gamma$ and $s_0$ and $r_0$ to be in $\Gamma_0$.\\
For $0 \leq i \leq n$ define $D_i(s_0,s) = t_0^{s_0} t^s t_i \frac{d}{dt_i}$ 
which acts on $A(m_0,m)$ as 
derivations. It is well known that $D(m_0,m)$ has the following basis
$$
\{D_i(s_0,s)| 0 \leq i \leq n, s_o \in \Gamma_0, s \in \Gamma \}.
$$
Let $d_i= t_i \frac{d}{dt_i}$ and it is easy to see that 

\item [{(2.1.1)}] $[t_0^{s_0} t^s d_a, t_0^{r_0} t^r d_b]= r_a t_0^{r_0+s_0} 
t^{r+s} d_b-s_b t_0^{r_0+s_0} t^{r+s} d_a .$

\item [{(2.2)}] $D(m_0,m)$ acts on $Z(m_0,m)$ in the following way 

\item [{(2.2.1)}] $t_0^{s_0} t^s d_a. (t_0^{r_0} t^r K_b)= r_a t_0^{r_0+s_0} 
t^{r+s} K_b + \delta_{ab} \displaystyle{\sum_{p=0}^n} s_p t_0^{r_0+s_0} t^{r+s} 
K_p .$

\item [{(2.3)}] It is known that $D(m_0,m)$ admits two non-trivial 
2-cocycles with values in $Z(m_0,m)$. See 
$[BB]$ for details
$$
\begin{array}{lll}
\varphi_1(t_0^{r_0} t^r d_a, t_0^{s_0} t^s d_b)&= -s_a r_b 
\displaystyle{\sum_{p=0}^n} r_p t_0^{r_0+s_0} t^{r+s} K_p,\\
\varphi_2(t_0^{r_0} t^r d_a, t_0^{s_0} t^s d_b)&= r_a s_b 
\displaystyle{\sum_{p=0}^n} r_p t_0^{r_0+s_0} t^{r+s} K_p.
\end{array}
$$

\item [{(2.4)}] Let $\varphi$ be arbitrary linear combinations of 
$\varphi_1$ and $\varphi_2$. Then there is a corresponding Lie algebra

\item [{(2.4.1)}] $\tau = L(\gg, \sigma) \oplus Z(m_0,m) \oplus D(m_0,m)$.\\
The Lie brackets  are defined in the following way in addition to 1.8.1 and 
1.8.2.

\item [{(2.4.2)}] $[t_0^{r_0} t^r d_a, X(k_0,k)] =k_a X(k_0+r_0, k+r),$

\item [{(2.4.3)}] $[t_0^{r_0} t^r d_a, t^{s_0} t^s K_b]= s_a t_0^{r_0+s_0} 
t^{r+s} K_b + \delta_{ab} \displaystyle{\sum_{p=0}^n} r_p t_0^{r_0+s_0} t^{r+s} 
K_p,$

\item [{(2.4.4)}] $[t_0^{r_0} t^r d_a, t_0^{s_0} t^s d_b] = s_a t_0^{r_0+s_0} 
t^{r+s} d_b-r_b  t_0^{r_0+s_0} t^{r+s} d_a + \varphi(t_0^{r_0} t^r d_a, t^{s_0} 
t^s d_b)$, \\
where $r,s \in \Gamma, r_0,s_0 \in \Gamma_0, X \in \gg (\overline{k}_0, 
\overline{k}).$
\end{enumerate}

\section{Assumptions and automorphisms.}

In this section we will make some assumptions on $L(\gg,\sigma)$ which will 
hold throughout this paper. We will also define a class of automorphisms on
$\tau$. 

\begin{enumerate}
\item [{(3.1)}] Assumptions

\item [{(3.1.1)}] $\gg(\overline{\circ}, \overline{\circ})$ is simple Lie 
algebra. 

\item [{(3.1.2)}] We can choose Cartan subalgebra $\hh(\circ)$ and $\hh$ for 
$\gg(\overline{\circ}, \overline{\circ})$ 
and $\gg$ such that $\hh(\circ) \subseteq \hh$.

\item [{(3.1.3)}] It is known that $\Delta^\times_0 = 
\Delta(\gg(\overline{\circ}, \overline{\circ}), \hh(\circ)) \backslash 
\{0\}$ is an irreducible reduced finite root system and has atmost two root lengths. Let $\Delta^\times_{0, sh}$ be the set of non-zero short roots. 
Define\\
$
\Delta^\times_{0,en}= \begin{cases} 
\Delta^\times_0 \cup 2 \Delta^\times_{0, sh} \ \mbox{if} \  \Delta^\times_0 \ 
\mbox{is of  \ type}\  B_l\\
\Delta^\times_0 \  \mbox{ otherwise}\\
\end{cases}
$                                              

$
\Delta_{0,en}= \Delta^\times_{0,en} \cup \{0\}.\\
$
We assume that $\Delta(\gg, \hh(0)) = \Delta_{0,en}$.

\item [{(3.2)}] {\bf Remark} These assumptions are true for any Lie Torus. 
See Proposition 3.2.5 of $[ABFP].$ 

It should be mentioned that Lie Torus are very important class of Lie 
algebras and give rise to almost all Extended Affine Lie algebras. See 
$[ABFP]$ and references there in. 

\item [{(3.3)}] Change of co-ordinates\\
In this subsection we will define a class of automorphisms for the Lie 
algebra
$$
\tau(1,{\bf 1}) = \gg \otimes A \oplus \Omega_A/d_A \oplus D(1,1).
$$

It is standard fact that $GL(n+1, \ZZ)$ acts on $\ZZ^{n+1}$ and we denote the
action by dot. Let 
$B=(b_{ij}) \in GL(n+1, \ZZ)$, then define automorphism, again denote by 
$B$ on $\tau(1,{\bf 1})$.\\
Let $t_0^{k_0} t^k= t (k_0,k)$, then define
$$
\begin{array}{lll}
B \cdot x \otimes t{(k_0,k)} &= x \otimes t^{B \cdot (k_0,k)},\\
B \cdot t{(k_0,k)} K_j &= \displaystyle{\sum_{p=0}^n} b_{pj} t^{B \cdot (k_0,k)} 
K_p,\\
B \cdot t{(k_0,k)} d_j &= \displaystyle{\sum_{p=0}^n} c_{pj} t^{B \cdot (k_0,k)} 
d_p,\\
\end{array}
$$
where $B^{-1}= (c_{pj})$.\

This  is what we call change of co-ordinates. We will use this change of 
co-ordinates without any mention 
and just say "upto choice of co-ordinates"
\end{enumerate}

\section{Root space decomposition and integrable modules 
for $\tau$.}

\begin{enumerate}
\item [{(4.1)}]  First note the center of $\tau$ is spanned by $K_0,K_1, 
\cdots, K_n$.\\ 
Let $H= \hh(0) \oplus \sum \CC K_i \oplus \sum \CC d_i$ which is an abelian 
Lie algebra of $\tau$ and plays the role of Cartan subalgebra.\\
Define $\delta_i, \Lambda_i \in H^* (0 \leq i \leq n)$ be such that 

\item [{(4.1.1)}]
$$
\begin{array}{lll}
\Lambda_i(\hh(0)) &= 0, \Lambda_i(K_j) = \delta_{ij}, \Lambda_i(d_j)=0,\\
\delta_i(\hh(0)) &= 0, \delta_i(K_j) = 0, \delta_i(d_j)= \delta_{ij}.
\end{array}\\
$$

Let $ \delta_k = \displaystyle{\sum_{i=1}^n} k_i \delta_i \ for\ k \in 
\ZZ^n$.

\item [{(4.1.2)}] Let $\gg(\overline{k}_0, \overline{k}, \alpha) = \{x \in 
\gg(\overline{k}_0, \overline{k}) | [h,x]= \alpha (h) x \ \mbox{for all}\ h
\in \hh(0)\}$\\
then $\tau$ has a root space decomposition.

\item [{(4.1.3)}] $\tau= \displaystyle{\bigoplus_{\beta \in \Delta}} \tau_\beta$ 

where $\Delta \subseteq \{\alpha + k_0 \delta_0 +\delta_k, \alpha \in 
\Delta_{0, e_n}, k_0 \in \ZZ, k \in \ZZ^n \}$.
$$
\begin{array}{lll}
\tau_{\alpha+k_0 \delta_0 +\delta_k} &= \gg(\overline{k}_0, \overline{k}, 
\alpha) \otimes t^{k_0}_0 t^k \ \mbox{for}\ \alpha \neq 0,\\
\tau_{k_0 \delta_0 +\delta_k} &= \gg(\overline{k}_0, \overline{k}, 0) \otimes
 t_0^{k_0} t^k \oplus
\displaystyle{\bigoplus^n_{i=0}} \CC t^{k_0}_0 t^k K_i \oplus 
\displaystyle{\bigoplus^n_{i=0}} \CC t^{k_0}_0 t^k d_i.
\end{array}\\
$$

Notice that $\tau_0=H$.

\item [{(4.1.4)}] Now we will define a non-degenerate bilinear form on $H^*$. 
For $\alpha \in \hh(0)^*$ extended 
$\alpha$ to $H$ by $\alpha(K_i)= \alpha(d_i) = 0, 0 \leq i \leq n$.\\
Let $(\hh(0), K_i)=0=(\hh(0),d_i)$,\\
$(\delta_k+\delta_{k_0}, \delta_l+\delta_{l_0}) =0= (\Lambda_k,\Lambda_p)$,\\
$(\delta_i, \Lambda_j) = \delta_{ij}$. The form on $\hh(0)$ is the 
restriction of the form $(,)$ on $\gg$.

\item [{(4.1.5)}] For $\gamma = \alpha+k_0 \delta_0+\delta_k$ is called real 
root if $\alpha \neq 0$ which is 
equivalent to $(\gamma,\gamma) \neq 0$. Denote $\Delta^{re}$ be the set of 
real roots. For 

$\alpha \in \Delta_{0,en}$, denote $\alpha^\vee$ the co-root of $\alpha$.\\
Define $\gamma^\vee= \alpha^\vee+ \frac{2}{(\alpha,\alpha)}
\displaystyle{\sum^n_{i=0}} k_i K_i$ for $\gamma$ real.\\
Then $\gamma(\gamma^\vee) = \alpha (\alpha^\vee)=2$.\\
For $\gamma$ real root, define reflection on $H^*$ by 
$$
r_\gamma(\lambda) = \lambda - \lambda(\gamma^\vee)\gamma, \gamma \in H^*.
$$

Let $W$ be the Weyl group genarated by $r_\gamma, \gamma \in \Delta^{re}$.

\item [{(4.2)}] A module $V$ of $\tau$ is called integrable if 

\item [{(4.2.1)}] $V= \displaystyle{\bigoplus_{\lambda \in H^*}} V_\lambda, 
V_\lambda =\{v \in V| h v=\lambda (h)v, h \in H\}, \dim V_\lambda < \infty$,

\item [{(4.2.2)}] $\gg(\overline{k}_0, \overline{k}, \alpha) \otimes 
t^{k_0}_0 t^k$ acts locally nilpotently on $V$
for $\alpha \neq 0$.\\
Let $P(V) = \{\gamma \in H^*| V_\gamma \neq 0\}$.\\
The following Lemma is very standard.

\item [{(4.3)}] {\bf Lemma} Suppose $V$ is an irreducible integrable module 
for $\tau$. Then 

\item [{(4.3.1)}] $P(V)$ is $W$- invariant. 

\item [{(4.3.2)}] $\dim V_\gamma = \dim V_{w\gamma}$ for all $w \in W$.

\item [{(4.3.3)}] For $\alpha \in \Delta^{re}, \lambda \in P(V)$ we have 
$\lambda(\alpha^\vee) \in \ZZ.$

\item [{(4.3.4)}] For $\alpha \in \Delta^{re}, \lambda \in P(V)$. If $\lambda 
(\alpha^\vee) > 0$ then 
$\lambda - \alpha \in P(V)$.

\item [{(4.3.5)}] For $\lambda \in P(V), \lambda (K_i)$ is a constant integer.

\end{enumerate}

The purpose of this paper is to classify irreducible integrable modules for 
$\tau$ with non-zero central action.

For an irreducible integrable module with non zero central charge, we can 
assume that $K_0$ acts as $C_0 > 0$ and $K_i (i \neq 0)$ acts trivially upto
a choice of co-ordinates. See (3.3).

\section{Existence of highest weight}

Throughout the rest of the paper we fix an irreducible integrable module 
for $\tau$ with $K_0$ acting as 
$C_0 >0$ and $K_i (i \neq 0)$ acts trivially.
Notice that for any $\lambda \in  P(V), \lambda(K_i)= C_i=0$ for $1 \leq i 
\leq n$ and $\lambda(K_0)=C_0$. 
Given a $\lambda \in H^*$ let $\lambda^\prime$ denote the restriction to 
$\hh(0)$. Given a $\lambda^\prime$ in 
$\hh^*(0)$, extend to $H$ by $\lambda^\prime(K_i)=\lambda(d_i)=0.$ Then 
$\lambda$ can be uniquely written as 
\begin{enumerate}

\item [{(5.1)}] $\lambda= \lambda^\prime + \displaystyle{\sum^n_{i=0}} 
\lambda (d_i) \delta_i+ \displaystyle{\sum^n_{i=0}} \lambda(K_i) \Lambda_i$.\\
For $\lambda \in P(V)$\\
$\lambda= \lambda^\prime + \lambda(d_0) \delta_0+ \displaystyle{\sum^n_{i=1}} 
\lambda(d_i) \delta_i + \lambda(K_0) \Lambda_0$.\\
Put $\overline{\lambda}= \lambda^\prime + \lambda(d_0) \delta_0+ \lambda(K_0)
 \Lambda_0$,\\
so that $\lambda= \overline{\lambda}+ \displaystyle{\sum^n_{i=1}} 
\lambda(d_i) \delta_i$.\
Let $\alpha_0= -\beta_0+\delta_0$ where $\beta_0$ is maximal root in 
$\Delta_{0,en}$. Note that $\alpha_0$ may not be root of $\tau$.\\
Let $\alpha_1, \alpha_2, \cdots, \alpha_p$ be a set of simple roots for 
$\Delta(g(\overline{\circ}, \overline{\circ}), \hh(0))$ and let $Q^+= 
\displaystyle{\bigoplus^p_{i=0}} \NN \alpha_i$.
Define an ordering on $H^*, \lambda \leq \mu$ for $\lambda, \mu \in H^*$, 
if $\mu-\lambda \in Q^+$. Notice that in this case 
$\mu(d_i)= \lambda(d_i)$ for $1 \leq i \leq n$.

\item [{(5.2)}] {\bf Proposition} Assumptions as above. Given a $\mu \in 
P(V)$ there exists a $\lambda \in P(V)$ 
such that $\mu \leq \lambda$ and $\lambda + \alpha \notin P(V)$ for all 
$\alpha > 0$, where $\lambda$ is dominant integral.\\
{\bf Proof} The proof of the Proposition follows from the results of 
Section 6 of [E4]. We will briefly sketch the proof. Suppose the Proposition
is false. Then the conclusion of Proposition 6.5 of  [E4] holds and it will 
lead to contradiction as explained in the proof of Theorem 6.1 of [E4]. Also
 note that the proof in [E4] is worked out only for $n=0$ but holds for any 
$n$. Here we need to observe in the construction of $\lambda_i$'s in 
Proposition 6.5 of  [E4], the $\delta_k, k \in \ZZ^n$ does not appear. In 
Proposition (5.2) $\lambda$ is clearly dominant.

\item [{(5.3)}] {\bf Theorem} Let $V$ be an irreducible integrable module 
for $\tau$ with $K_0$ acts as 
$C_0 > 0$ and $K_i$ acts trivially. Then there exists $\gamma \in P(V)$ 
such that $\gamma + \beta + \delta_k \notin P(V)$ 
for any $\beta >0$ and for any $k \in \ZZ^n$.\\
{\bf Proof} Suppose the theorem is false. First some word about notation. 
The $\delta^\prime$s that occur below 
are  always linear combinations of $\delta_1, \cdots, \delta_n$. Fix 
$\lambda_1 \in P(V)$ and let 
$\delta(g)= \displaystyle{\sum_{i=1}^n} \lambda(d_i) \delta_i$. Note that 
$(\lambda_1, \delta_i)=0$ for 
$1 \leq i \leq n$ which follows from unique expression 5.1. Choose $\mu_1$ 
as in the Proposition (5.2). Then 
$\lambda_1 \leq \mu_1$ and $\mu_1+ \alpha \notin P(V)$ for $\alpha >0$. By 
assumptions there exists a root 
$\beta_1 + \delta(1)$ with $\beta_1 > 0$ and $\delta(1) \neq 0$ such that
$$
\lambda_2= \mu_1 +\beta_1+\delta(1) \in P(V).
$$
Notice that $\lambda_2= \overline{\lambda}_2+ \delta(g) + \delta(1)$ and 
$\overline{\lambda}_1 \leq \overline{\mu}_1  < \overline{\lambda}_2$.\\
By repeating  the argument infinitely many times we get dominant integral 
weights $\mu_d \in P(V)$ and roots $\beta_d + \delta(d)$ such that 
$\beta_d > 0$ and $\delta(d) \neq 0$ such that

\item [{(5.3.1)}] $\lambda_{d+1} = \mu_d + \beta_{d} + \delta(d) \in P(V)$. 

\item [{(5.3.2)}] $\mu_{d} + \beta \notin P(V)$ for all positive $\beta$.

\item [{(5.3.3)}]  $\lambda_{d} \leq \mu_{d}$, 

\item [{(5.3.4)}] $\overline{\lambda}_{d }    \leq  \overline{\mu}_{d}
 <  \overline{\lambda}_{d+1} \leq  \overline{\mu}_{d+1}$,\\
It is easy to see that
\item [{(5.3.5)}]$\mu_{d} = \overline{\mu}_{d} + \delta(g) + \delta^\prime (d)$, 
where $\delta^\prime (d)= \displaystyle{\sum_{k=1}^d} \delta(k)$.\\
We also have 
\item [{(5.3.6)}] $\overline{\mu}_{d_1} < \overline{\mu}_{d_2}$ for $d_1 < d_2$.

\item [{(5.3.7)}] 
$(\mu_{d} + \beta_{d} , \beta_{d}) >0 $, as $\mu_{d}$ is dominant integral.

\item [{(5.3.8)}] {\bf Claim} $\mu_{d_1} < \mu_{d_2}- \delta^\prime (d_2) + 
\delta^\prime (d_1)$ for $d_1 < d_2$. 
\end{enumerate}

To see the claim note that from (5.3.5) it follows that 
$$
\mu_{d_1}- \delta(g)- \delta^\prime (d_1) < \mu_{d_2} - \delta(g) - 
\delta^\prime (d_2).
$$

The inequality still holds if we add $\delta(g) + \delta^\prime(d_1)$ both 
sides. The claim follows.\

Now there exists $d_1 < d_2$ such that $\delta^\prime (d_2) - \delta^\prime 
(d_1) = \delta_s$ for some $s \in \Gamma$. We know that $\beta_{d_2} + 
\delta(d_2) \in \Delta$. It follows from the definition of $\tau$ that 
$$
\beta_{d_2} + \delta(d_2) + \delta^\prime (d_2) - \delta^\prime (d_1) \in 
\Delta.
$$

Now from  (5.1) and (5.3.7) it follows that 
$$
(\mu_{d_2} + \beta_{d_2} + \delta(d_2), \beta_{d_2}+ \delta(d_2) + 
\delta^\prime (d_2) - \delta^\prime (d_1) >0.
$$

From (4.3.4) it follows that 

$$
\mu_{d_2} + \delta^\prime (d_1)-\delta^\prime (d_2) \in P(V).
$$
But by claim it follows that 
$$
\mu_{d_1} < \mu_{d_2} + \delta^\prime (d_1) - \delta^\prime (d_2),
$$
which is a contradiction to (5.3.2). This proves the Theorem.

\section{Triangle decomposition}

\begin{enumerate}

\item [{(6.1)}] We will now define triangular decomposition for $\tau$. Let 
$Z= \Omega A/dA$.\\
Let 
$$
\begin{array}{lll}
L^+(\gg, {\bf \sigma}) &= \displaystyle{\bigoplus_{\alpha+ k_0 \delta_0 >0}}
\gg(\overline{k}_0, \overline{k}, \alpha)\otimes  t_0^{k_0} t^k, k \in 
\ZZ^n;\\
L^-(\gg, {\bf \sigma}) &= \displaystyle{\bigoplus_{\alpha+ k_0 \delta_0 <0}}
\gg(\overline{k}_0, \overline{k}, \alpha)\otimes  t_0^{k_0} t^k, k \in 
\ZZ^n;\\
L^0(\gg, {\bf \sigma}) &= \displaystyle{\bigoplus_{k \in \ZZ^n}} 
\gg(\overline{0}, \overline{k}, 0) t^k;\\
D^+(m_0,m) &= \displaystyle\bigoplus_{\substack{0 \leq i \leq n\\ s_0 > 0}} 
\CC t_0^{s_0} t^s d_i, s \in \Gamma;\\
D^-(m_0,m) &= \displaystyle\bigoplus_{\substack{0 \leq i \leq n\\ s_0 < 0}} 
\CC t_0^{s_0} t^s d_i, s \in \Gamma;\\
D^0(m_0,m) &= \displaystyle\bigoplus_{0 \leq i \leq n} \CC t^s d_i, s \in 
\Gamma;\\
Z^+ &= \displaystyle\bigoplus_{\substack{0 \leq i \leq n\\ s_0 > 0}} \CC 
t_0^{s_0} t^s K_i, s \in \Gamma;\\
Z^- &= \displaystyle\bigoplus_{\substack{0 \leq i \leq n\\ s_0 < 0}} \CC 
t_0^{s_0} t^s K_i, s \in \Gamma;\\
Z^0 &= \displaystyle\bigoplus_{0 \leq i \leq n} \CC t^s K_i, s \in \Gamma;\\
\tau^+ &= L^+(\gg, {\bf \sigma}) \oplus Z^+ \oplus D^+ (m_0,m);\\
\tau^- &= L^-(\gg, {\bf \sigma}) \oplus Z^- \oplus D^- (m_0,m);\\
\tau^0 &= L^0(\gg, {\bf \sigma}) \oplus Z^0 \oplus D^0 (m_0,m).\\
\end{array}
$$

Then  clearly $\tau = \tau^- \oplus \tau^0 \oplus \tau^+$ is a triangular 
decomposition.\\
Let $T= \{ v \in V| \tau^+ v = 0\} \neq 0$ by Theorem 5.3

\item [{(6.2)}] {\bf Proposition:} $T$ is a $\tau^0$- module and in fact 
irreducible as $\tau^0$- module. 
Further $V= U(\tau^-) T.$\\
{\bf Proof} It is easy to check that  $[\tau^0, \tau^+] \subset \tau^+$. 
From this  it follows that  $T$ is a $\tau^0$- module. Now from  PBW 
theorem, we have $U(\tau) =  U(\tau^-) U( \tau^0)U(\tau^+)$.  ( Here $U$ 
denotes  the universal enveloping algebra.) Using this and the fact that  
$V$ is $\tau$-irreducible, it follows that $T$ is  $\tau^0$-irreducible and 
$V =  U(\tau^-) T$. Recall that 
$\{d_1, \cdots, d_n\} \subseteq D^0 (m_0,m)$ and hence $T$ is $\ZZ^n$- 
graded.\\
Let $T_k= \{v \in T| d_i v = (\lambda(d_i)+ k_i)v, 1 \leq i \leq n \}$ 

where $\lambda$ is a fixed weight in $P(V)$. We will now  record some 
result on $T$ which can be proved similarly as in $[JM], [EJ]$ and $[FJ]$.

\item [{(6.3)}] {\bf Lemma}

\item [{(6.3.1)}] For $v \in T \backslash \{0\}, t^s K_0 v \neq 0$ for all 
$s \in \Gamma$. 

\item [{(6.3.2)}] $\dim T_k = \dim T_{k+s} = d_k \ \forall \ s \in \Gamma$.

\item [{(6.3.3)}] Let $v_1 (k), v_2(k), \cdots, v_{d_k} (k)$ be a  basis of 
$T_k$. Let 
$v_i(s+k) = \frac{1}{C_0} t^s K_0 v_i(k)$, then $\{v_1(s+k), \cdots, v_{d_k} 
(s+k)\}$ is a basis of $T_{k+s}$.

\item [{(6.3.4)}] $\frac{1}{C_0} t^s K_0 (v_1(k+r), \cdots, v_{d_k} (k+r)) = 
(v_1(k+r+s), \cdots, v_{d_k} (k+r+s))$ for all $r, s \in \Gamma$.

\item [{(6.3.5)}] $h \otimes t^s(v_1(k+r), \cdots, v_{d_k} (k+r))= \lambda (h)
(v_1(k+r+s), \cdots, v_{d_k} (k+r+s))$ for $h \in \hh(0)$ and all $r, s \in 
\Gamma$, where $\lambda$ is a fixed weight of $P(V)$.

\item [{(6.3.6)}] $t^s d_0 (v_1(k+r), \cdots, v_{d_k} (k+r))= \lambda (d_0) 
(v_1(k+r+s), \cdots, v_{d_k} (k+r+s))$ for all $s, r \in \Gamma$ and for a 
fixed $\lambda \in P(V)$.

\item [{(6.3.7)}] $t^r K_p \cdot T=0 \ \  1 \leq p \leq n, r \in \Gamma$.

\item [{(6.3.8)}] $t^r K_0 \cdot t^s K_0 \cdot v= C_0 t^{r+s} K_0 v \ \ \forall
 v \in T$ and $r, s \in \Gamma$.
\end{enumerate}

\section{More notation and co-finite ideals}

\begin{enumerate}

\item [{(7.1)}]Let $\Der A(m)$ be the derivation algebra of $A(m)$. Let 
$e_1, \cdots, e_n$ be the standard basis 
of $\CC^n$ and let $u= \sum u_i e_i \in \CC^n$. Let $D(u,r)= 
\displaystyle{\sum_{i=1}^n} u_i t^r d_i, r \in \Gamma$.

\item [{(7.2)}] From the earlier section $T$ can be identified with $V^1 
\otimes A(m)$ where $V^1$ can be identified with 
$$
\substack{\bigoplus {\bf T}_k\\ 0 \leq k_i < m_i\\ 1 \leq i \leq n}
$$

\item [{(7.3)}] Now $D^0(m_0,m)$ is spanned by $t^r d_i, r \in \Gamma, 0 
\leq i \leq n$. Thus $D^0(m_0,m)$ can be identified with $Der A(m) \oplus 
\displaystyle{\sum_{r \in \Gamma}} \CC t^r d_0$\\

$Z^0$ can be identified with $\displaystyle{\sum_{r \in \Gamma}} \CC t^r K_0$ 
as the rest of the space acts trivially on $T$. Thus $V^1 \otimes A(m)$ is 
an irreducible module for 
$$
L= L^0(\gg, {\bf \sigma}) \oplus Der A(m) \oplus 
\displaystyle{\sum_{r \in \Gamma}} \CC t^r d_0 \oplus A(m),
$$
where we identify $\displaystyle{\sum_{r \in \Gamma}} \CC t^r K_0$ by $A(m)$ 
and $\frac{1}{C_0} t^r K_0$ goes to $t^r$ which is well defined by 6.3.8.

\item [{(7.4)}] We note the following 
$$
\begin{array}{lll}
t^r \cdot v \otimes t^s &= v \otimes t^{r+s},\\
t^r d_0 \cdot v \otimes t^s &= \lambda(d_0) v \otimes t^{r+s} \ \mbox{for} 
\ r, s \in \Gamma, v \in V^1.\\
\end{array}
$$

\item [{(7.5)}] Let $\stackrel{\circ}{\gg}= \{ X \in \gg | \sigma_0 X = X,
[h,X]=0, h \in \hh(0)\}$\\
the following is easily checked. 

\item [{(7.5.1)}] $\sigma_i (\stackrel{\circ}{\gg}) \subseteq 
\stackrel{\circ}{\gg} $ for $1 \leq i \leq n$.

\item [{(7.5.2)}] $\stackrel{\circ}{\gg} = 
\displaystyle{\bigoplus_{\overline{k} \in \Lambda}} \stackrel{\circ}
{\gg}_{\overline{k}}$ is a natural\\ $\Lambda$- grading where 
$\stackrel{\circ}{\gg}_{\overline{k}}= \{ X \in \stackrel{\circ}{\gg}| 
\sigma_i X = \xi_i^{k_i} X, 1 \leq i \leq n \}$\\
The corresponding multiloop algebra is denoted by \\
$L(\stackrel{\circ}{\gg}, \sigma) = \displaystyle{\bigoplus_{k \in \Lambda}}
\stackrel{\circ}{\gg}_{\overline{k}} \otimes t^k$\\
It is clear that $L^0(\gg, \sigma) = L(\stackrel{\circ}{\gg}, \sigma).$\\
When we say $X(k)= X \otimes t^k \in L(\stackrel{\circ}{\gg}, \sigma) $ we 
always mean $X \in \stackrel{\circ}{\gg}_{\overline{k}}$. \\
Thus $L \cong L(\stackrel{\circ}{\gg}, \sigma) \oplus Der A(m) \oplus A(m) 
\oplus \displaystyle{\sum_{r \in \Gamma}} \CC t^r d_0$.

\item [{(7.6)}] The brackets in $L$ are given as follows :

\item [{(7.6.1)}] $[X(k), Y(l)]= [X,Y] (k+l),$

\item [{(7.6.2)}] $[D(u,r), D(v,s)] = D(w, r+s)$ where $w= (u,s) v - (v,r) u,$

\item [{(7.6.3)}] $[D(u,r), t^s] = (u,s) t^{r+s},$

\item [{(7.6.4)}] $[D(u,r), X(k)]= (u,k) X (k+r),$

\item [{(7.6.5)}] $[D(u,r), t^s d_0]= (u,s) t^{r+s} d_0$.\

Now we would like to classify the irreducible $L$- module $V^1 \otimes 
A(m)$. We need some preparation for that.

\item [{(7.7)}] For $k \in \ZZ^n, X \in \stackrel{\circ}{\gg}_{\overline{k}}$ 
and $r_1, r_2, \cdots, r_d \in \Gamma$, define

$X(k, r_1, \cdots, r_d)= X(k) - \sum X(k+r_i) + \displaystyle{\sum_{i < j}}
X(k+r_i+r_j) + \cdots (-1)^d X (k+r_1 + r_2 + \cdots + r_d)$\\
Let $F_d$ be the subspace of $L(\stackrel{\circ}{\gg}, \sigma) $ spanned by
$X(k, r_1, \cdots, r_d)$. It is easily checked that $F_d$ is an ideal in 
$L(\stackrel{\circ}{\gg}, \sigma)$

\item [{(7.8)}] {\bf Lemma}

\item [{(7.8.1)}] $F_d \subseteq F_{d-1}$

\item [{(7.8.2)}] $[F_d, F_d] \subseteq F_{d+1}$\\
{\bf Proof-} Note that $X(k, r_1, \cdots, r_d) = X(k, r_1, \cdots, r_{d-1}) -
X(k+ r_d, r_1, \cdots, r_{d-1})$\\
which proves (7.8.1).\\
Now consider for $l, k \in \ZZ^n, r_1, \cdots, r_d, s \in \Gamma, X \in 
\stackrel{\circ}{\gg}_{\overline{k}}, Y \in 
\stackrel{\circ}{\gg}_{\overline{l}}$.\\
$[X(k, r_1, \cdots, r_d), Y(l) - Y(l+s)]$\\
$=[X,Y] (k+l, r_1, \cdots, r_d)$\\
$-[X,Y] (k+l+s, r_1, \cdots, r_d)$\\
$=[X,Y] (k+l, r_1, r_2, \cdots, r_d,s) \in F_{d+1}$.\\
By the above note we see that $F_d$ is spanned by vectors $X(k) - X(k+s)$. 
Thus from above (7.8.2) follows.

\item [{(7.9)}] In this subsection we recall some facts from $[E2]$ on $Der 
A(m)$. Let $I(u,r)= D(u,r)- D(u,0), u \in \CC^n, r \in \Gamma$ It is easy to
check.

\item [{(7.9.1)}]

$
[I(u,r), I(v,s)] = (v,r) I(u,r) -(u,s) I(v,s) + I(w,s+r)\\
$
where $w= (u,s) v-(v,r)u$

Let  $I$ be the space spanned by $I(u,r), u \in \CC^n, r \in \Gamma$ which 
can be seen as subalgebra of $Der A(m)$.

\item [{(7.10)}] For $d \geq 1, u \in \CC^n, s_1, \cdots, s_d, r \in 
\Gamma, $\\
Let  $I_d(u,r,s_1,s_2, \cdots, s_d)= I(u,r)- \displaystyle{\sum_i}
I(u,r+s_i)+ \displaystyle{\sum_{i<j}} I(u,r+s_i+s_j) + \cdots + (-1)^d 
I(u,r+s_1+s_2 +  \cdots + s_d)$\\
Let $I_d$ be the space spanned by  $I_d(u,r,s_1,s_2, \cdots, s_d)$ for \\
$u \in \CC^n, r,s_1, \cdots s_d \in \Gamma$.\\
The following is proved in $[E2]$

\item [{(7.11)}] {\bf Lemma} 

\item [{(7.11.1)}]$I_d$ is a co-finite ideal in $I$.

\item [{(7.11.2)}] Any co-finite ideal of $I$ contains $I_d$ for large $d$.

\item [{(7.11.3)}] $I_1=I$ and $I/I_2 \cong gl_n$.
\end{enumerate}
\section{Finite dimensional modules}
\begin{enumerate}

\item [{(8.1)}]  Let $W$ be the subspace of $V^1 \otimes A(m)$ spanned by 
vectors of the form
$t^r. v(s)-v(s)$ for $r,s \in \Gamma$ and $v \in V^1$.\\
Let $\stackrel{\sim}{L}= I \ltimes L(\stackrel{\circ}{\gg}, \sigma) $

\item [{(8.2)}] {\bf Lemma-} $W$ is an $\stackrel{\sim}{L} \oplus A(m) 
\oplus \displaystyle{\sum_{r \in \Gamma}} \CC t^r d_0$ module.\\
{\bf Proof-} It is easy to check using the following
$$
\begin{array}{lll}
[D(u,r)- D(u,0), X(k)] &= (u,k) (X(k+r)-X(k)),\\
{[D(u,r)]- D(u,0), t^s]}&= 0,\\
{[L(\stackrel{\circ}{\gg}, \sigma), A(m)]}& = 0,\\
{[t^r d_0, A(m)]} &=0.
\end{array}
$$

\end{enumerate}

Let $\stackrel{\sim}{V}= (V^{1} \otimes A(m))/W$ which is an 
$\stackrel{\sim}{L}$- module. Notice that 
$A(m) \oplus \displaystyle{\sum_{r \in \Gamma}} \CC t^r d_0$ acts as  scalars 
on $\stackrel{\sim}{V}$ and hence 
we ignore them.
We would like to prove that $\stackrel{\sim}{V}$ is completely reducible 
$\stackrel{\sim}{L}$-module.

\begin{enumerate}

\item [{(8.3)}] Recall that the Lie brackets in $\stackrel{\sim}{L}$ are 
given by 

\item [{(8.3.1)}] $[I(v,s), I(u,r)]= (u,s) I(v,s) - (v,r) I(u,r) + I(w, 
r+s)$,\\
where $I(u,r)= D(u,r)-D(u,0)$ and $w=(v,r)u-(u,s) v$

\item [{(8.3.2)}] $[I(v,s), X(k)] = (v,k) (X(s+k)-X(k))$,

\item [{(8.3.3)}] $[X(k), Y(l)] = [X,Y](k+l)$,
\end{enumerate}

where $X \in \gg_{\overline{k}}, Y \in \gg_{\overline{l}}, k,l \in \ZZ^n, r, s 
\in \Gamma$ and $u, v \in \CC^n$.

\begin{enumerate}

\item [{(8.4)}] Recall that we fixed $\lambda \in P(V).$ Let $\alpha_i= 
\lambda (d_i)$ 
Let $\alpha= \sum \alpha_i e_i \in \CC^n$ and let $V_1$ is an 
$\stackrel{\sim}{L}$- module. Then we will define 
$L$- module structure on $L(V_1) = V_1 \otimes A_n$.
$$
\begin{array}{lll}
X(k) \cdot v_1 \otimes t^l &=  (X(k) v_1) \otimes t^{l+k},\\
D(u,r) \cdot v_1 \otimes t^l &=  (I(u,r) v_1) \otimes t^{l+r} + (u,l+\alpha)
v_1 \otimes t^{l+r},\\
t^s v_1 \otimes t^l &=  v_1 \otimes t^{s+l},\\
t^r d_0 \cdot v_1 \otimes t^l &= \lambda (d_0) \cdot v_1 \otimes t^{l+r},
\end{array}
$$

\end{enumerate}
where $v_1 \in V_1, l, k \in \ZZ^n, r,s, \in \Gamma.$\\
We need to check the brackets in (7.6). We will first check (7.6.4).
Consider\\
$D(u,r) X(k) (v_1 \otimes t^l)= D(u,r) ((X(k) v_1 \otimes t^{l+k})$\\
$= I(u,r) \cdot X(k) v_1 \otimes t^{l+k+r} + (u, l+k+ \alpha) X(k) v_1 
\otimes t^{l+k+r}$.\\
Consider\\
$X(k) D(u,r) (v_1 \otimes t^l)= X(k) (I(u,r) v_1 \otimes t^{l+r} + 
(u,l+\alpha) v_1 \otimes t^{l+r})$\\
$= X(k) I(u,r) v_1 \otimes t^{l+k+r} + (u,l+\alpha) X(k) v_1 
\otimes t^{l+k+r}$.\\
Now we will use the fact $[I(u,r), X(k)] = (u,k) (X(k+r) - X(k))$.\\
So \  $(D(u,r) X(k) - X(k) D(u,r)) \cdot (v_1 \otimes t^l)$\\
$= (u,k) (X(k+r) -X(k)) v_1 \otimes t^{l+k+r} + (u,k) X(k) v_1 \otimes 
t^{l+k+r}$\\
$=(u,k) X(k+r) v_1 \otimes t^{l+k+r}$\\
$=(u,k) X(k+r) (v_1 \otimes t^l)$.\\
This proves (7.6.4).\\
We will now check (7.6.2).\\
Consider\\
$D(v,s) D(u,r) (v_1 \otimes t^l) = D(v,s) (I(u,r) v_1 \otimes t^{l+r} + 
(u, l+\alpha) v_1 \otimes t^{l+r}$)\\
$=I(v,s) I(u,r) v_1 \otimes t^{l+r+s} + (v, l+r+\alpha) I(u,r) v_1 \otimes 
t^{l+r+s}$\\
$+ (u, l+\alpha) I(v,s) v_1 \otimes t^{l+r+s} + (u, l+\alpha) (v, l+r+\alpha)
 v_1 \otimes t^{l+r+s}$.\\
Similarly we have\\
$D(u,r) D(v,s) (v_1 \otimes t^l)$\\
$=I(u,r) I(v,s) v_1 \otimes t^{l+r+s}$\\
$+ (u,l+s+ \alpha) I(v,s) v_1 \otimes t^{l+r+s}$\\
$+(v, l+\alpha) I(u,r) v_1 \otimes t^{l+r+s}$\\
$+(v,l+\alpha) (u, l+s+\alpha) v_1 \otimes t^{l+r+s}$.\\
Now we will use (8.3.1)\\
So\\
$
(D(v,s) D(u,r) - D(u,r) D(v,s)) v_1 \otimes t^l\\
=((u,s) I(v,s) - (v,r) I(u,r) +I(w,r+s)) v_1 \otimes t^{l+r+s}\\
+(v,r) I(u,r) v_1 \otimes t^{l+r+s}\\
-(u,s) I(v,s) v_1 \otimes t^{l+r+s}\\
+(w, l+\alpha) v_1 \otimes t^{l+r+s}\\
=D(w, r+s) (v_1 \otimes t^l)
$\\
This  proves (7.6.2)\\
The remaining brackets 7.6.1, 7.6.3, 7.6.5 are trivial to verify.

\begin{enumerate}

\item [{(8.5)}] Recall

\item [{(8.5.1)}] $T$ is an irreducible $L$-module from Proposition 6.2.

\item [{(8.5.2)}] $\stackrel{\sim}{V}$ is an $\stackrel{\sim}{L}$-module from 
Lemma 8.2.

\item [{(8.5.3)}] $L(\stackrel{\sim}{V})$ is an $L$-module from (8.4).

We will now establish that $T$ is contained in $L(\stackrel{\sim}{V})$ as 
$L$-modules.

For $v_k \in T_k$, let $\overline{v}_k$ be the image in $\stackrel{\sim}{V} 
\cong T/W$.\\
Let $\stackrel{\sim}{\varphi} : T \rightarrow L(\stackrel{\sim}{V})$\\
$\stackrel{\sim}{\varphi} (v_k) = \overline{v}_k \otimes t^k, k \in \ZZ^n$.

\item [{(8.6)}] {\bf Lemma} $\stackrel{\sim}{\varphi}$ is an $L$-module 
map.\\
{\bf Proof} Consider $\stackrel{\sim}{\varphi} (D(u,r) v_k) = 
\overline{D(u,r) v_k} \otimes t^{k+r}$.\\
Now
$$
\begin{array}{lll}
D(u,r) (\overline{v}_k \otimes t^k) &= \overline{(D(u,r) - D(u,0))v_k} 
\otimes t^{k+r}\\
&+ (u,k+\alpha) \overline{v}_k \otimes t^{k+r}\\
&= \overline{D(u,r) v_k} \otimes t^{k+r}.\\
\end{array}
$$

\end{enumerate}

Thus  we have verified $\stackrel{\sim}{\varphi} (D(u,r) v_k) = D(u,r) 
\stackrel{\sim}{\varphi} (v_k)$.\\

The rest of the relations are easy to verify.\
Cleary $\stackrel{\sim}{\varphi}$ is a non-zero map and since $T$ is an 
irreducible $L$-module we have 
$T \subseteq L(\stackrel{\sim}{V})$ as $L$-module.

\begin{enumerate}

\item [{(8.7)}] {\bf Theorem} $\stackrel{\sim}{V}$ is completely reducible
as  $\stackrel{\sim}{L}$- module and all components are isomorphic.\\
We will prove some results before proving Theorem (8.7).

\item [{(8.8)}] {\bf Proposition} Let $\mathfrak g$ be a Lie algebra. Let 
$V_1, V_2, \cdots ,V_n$ be mutually non-isomorphic irreducible  $\mathfrak
 g$-modules. Suppose $W$ is a non-zero $\mathfrak g$-submodule of 
$\displaystyle{\bigoplus_{i =1}^n V_i}$. Then there exists $S \subset \{1, 2,
\cdots, n \}$ such that
                  $$W = \displaystyle{\bigoplus_{i \in S} V_i}.$$
{\bf Proof}  Let $\stackrel{\sim}{\pi}_j : \displaystyle{\bigoplus_{i =1}^n 
V_i } \rightarrow V_j$ be the natural projection. Let ${\pi}_j = 
\stackrel{\sim}{\pi}_j| W$. Let
$S = \{ j | {\pi}_j \neq 0 \}$. Suppose $j \in S$, then $\pi_j(W) \neq 0 $.
Since $V_j$ is irreducible, it follows that $\pi_j(W) =  V_j$. Clearly
 $$W \subset \displaystyle{\bigoplus _{j \in S} V_j}.$$
{\bf Claim} $$W = \displaystyle{\bigoplus _{j \in S} V_j}.$$
Note that the claim proves the proposition. We prove the claim by induction 
on $n$. Certainly the claim is true for $n = 1$. Let $w \in W$ and write 
$w = v_{i_1} + v_{i_2} +...+ v_{i_k}$, where $0 \neq v_{i_j} \in V_{i_j}.$ Then 
we define $l(w) = k$. We will now prove the claim for $n=2$. Suppose $l(w) 
= 1$ for all $w \in W$, then clearly $W = V_1$ or $W = V_2$ and hence we are
done. Suppose
$l(w_0) = 2$ for some $w_0 \in W$. Then $w_0 = v_1 + v_2$ and $0 \neq v_i \in 
V_i$ for $i = 1, 2$. As noted earlier we have ${\pi}_1(W) = V_1$  and 
${\pi}_2(W) = V_2$. Suppose ker${\pi}_1$ = ker${\pi}_2 = 0$, then $W = 
{\pi}_i(W) = V_i$ for $i = 1, 2$. This proves that $ V_1 \cong  V_2$  as 
$\mathfrak g$-modules, which is a contradiction.
Now suppose ker${\pi}_1 = 0$. Then there exists $w \in$ ker${\pi}_1$ and 
$w = v^1 + v^2$, $ v^1 \in V_1$ and $ v^2 \in V_2$. But $0 = {\pi}_1(w) = 
v^1$. 
Hence $v^2 \in W$. This proves $V_2 \subset W$.
Recall that $w_0 = v_1 + v_2 \in W$. It now follows that $v_1 \in W$ and 
hence
$V_1 \subset W$. Thus $V_1 + V_2 = W$ and we are done. This completes the 
claim for $n=2$. 

Now assume the claim for $n-1$ and we prove it for $n$. Note that $n \geq 3$.
Suppose $S \subsetneqq \{1, 2, \cdots, n \}$. Then as noted earlier 
$W  \subsetneqq  \displaystyle{\bigoplus _{i \in S} V_i}$ and by induction the 
claim follows.
We can now assume that $S  = \{1, 2, \cdots, n \}$. Suppose ker${\pi}_i = 0$
for all $i$. Then $W = {\pi}_j(W)= V_j$ for all $j$. This proves that
$ V_j \cong  V_i$ for all $i$ and $j$, which is a contradiction to our 
assumption. We can now assume  that ker${\pi}_i \neq 0$ for some $i$. Let 
$w \in$  ker${\pi}_i$ and write $w = v_1 + v_2 + \cdots + v_n, v_j \in V_j$. 
Now $v_i = {\pi}_i(w) = 0$. This proves $l(w) < n$. Let $W_1$ be the 
submodule generated by
$w$. Then clearly $W_1 \subset \displaystyle{\bigoplus _{j \neq i} V_j}$. 
By induction it follows that there exists $T_1  \subsetneqq  \{1, 2, \cdots, 
n \}$
such that $W_1 = \displaystyle{\bigoplus _{j \in T_1} V_j} \subsetneqq W$. Let 
$T_2$ be the maximal subset of $\{1, 2, \cdots, n \}$ such that 
$\displaystyle{\bigoplus _{j \in T_2} V_j} \subset W$.
To prove the claim it is sufficient to prove $T_2 = \{1, 2, \cdots, n \}$. 
So suppose there exists $j \notin T_2$. Since ${\pi}_j(W) \neq 0$, there 
exists $w \in W$ such that ${\pi}_j(w) \neq 0$. This proves $ w \notin 
\displaystyle{\bigoplus _{k \in T_2} V_k}$. Now it is easy to find $w_2 \in 
W$ such that $w_2 =  v_{i_1} +  \cdots + v_{i_l}$ and $0 \neq v_{i_j} \in 
V_{i_j}$ and $\{i_1, i_2, \cdots, i_l \} \cap T_2 = \emptyset$. Now let $W_2$
be the submodule generated by $w_2$. Note that  $W_2 \subset 
\displaystyle{\bigoplus _{j = 1}^l V_{i_j}}$. 
Now by induction we can find $T_3$ such that $T_3 \cap T_2 = \emptyset$ and 
$W_2 =  \displaystyle{\bigoplus _{k \in T_3} V_k}$. Thus 
$\displaystyle{\bigoplus _{i \in T_2 \cup T_3} V_i} \subset W$ which 
contradicts the maximality of $T_2$. Thus $T_2 = \{1, 2, \cdots, n \}$ and 
proves the claim.
\item [{(8.9)}] {\bf Lemma}  $\stackrel{\sim}{V}$ is graded irreducible 
$\stackrel{\sim}{L}$-module.

{\bf Proof}  Recall
that $\stackrel{\sim}{V}$ and $\stackrel{\sim}{L}$ are $\Lambda$-graded. 
Further using the map  $\stackrel{\sim}{\varphi}$, $T 
\cong \stackrel{\sim}{\varphi}(T)$
as $L$-modules. We will first note that the $L$-module generated by 
$\overline{v}_k \otimes t^k, k \in \ZZ^n$ is equal to the $\stackrel{\sim}
{L} \oplus A(m)$-module generated by 
$\overline{v}_k \otimes t^k, k \in \ZZ^n$. This follows from the fact
$$D(u,r)( \overline{v}_k \otimes t^k)  =  \overline{I(u,r) v_k} \otimes 
t^{k+r} + (u,k+\alpha) \overline{v}_k \otimes t^{k+r}$$ and
$$t^r. {\overline{v}_k} \otimes t^k =  {\overline{v}_k}  \otimes t^{k+r}.$$
Also note that $t^r$ acts trivially on $\stackrel{\sim}{V}$. From this it is
easy to see that $\stackrel{\sim}{V}$ is graded irreducible.

We will now prove a decomposition theorem for $L(\stackrel{\sim}{V})$ as 
$L$-module. First we give some notation.
Recall that  $$\stackrel{\sim}{V} =  \displaystyle{\bigoplus _{{\overline p} 
\in \Lambda } \stackrel{\sim}{V}_{\overline p}}.$$
Let $p \in \ZZ^n$ and  ${\overline p} \in \Lambda$. Define
$$L(\stackrel{\sim}{V})(\overline p) = \{ {\overline{v}_k} \otimes t^{k+r+p},
{\overline{v}_k} \in \stackrel{\sim}{V}_{\overline k}, 
r \in \Gamma, k \in \ZZ^n \}.$$
Clearly $L(\stackrel{\sim}{V})(\overline p)$ is closed under $A(m)$ and 
$\displaystyle{\sum_{r \in \Gamma}} \CC t^r d_0$. \\
Consider for $X \in \stackrel{\circ}{\gg}_l$, 
$$X(l). ({\overline{v}_k} \otimes t^{k+r+p}) = \overline{X(l){v}_k} \otimes 
t^{k+l+r+p} \in L(\stackrel{\sim}{V})(\overline p).$$
$$D(u,s) ({\overline{v}_k} \otimes t^{k+r+p}) = \overline{I(u,s) v_k} \otimes 
t^{k+r+p+s} + (u, k+r+p +\alpha) {\overline{v}_k} \otimes t^{k+r+p+s}, s \in 
\Gamma.$$
We see that the above vector belongs to $L(\stackrel{\sim}{V})(\overline p)$.
Thus $L(\stackrel{\sim}{V})(\overline p)$ is an $L$-module.
Clearly $$L(\stackrel{\sim}{V}) = \displaystyle{\bigoplus _{{\overline p} 
\in \Lambda }L( \stackrel{\sim}{V})(\overline p)}$$
which is a finite sum of $L$-modules. We have seen already that 
$T \cong L(\stackrel{\sim}{V})(\overline 0)$ as $L$-modules and in particular
$L(\stackrel{\sim}{V})(\overline 0)$ is an irreducible $L$-module.

\item [{(8.10)}] {\bf Proposition} 
Each $L( \stackrel{\sim}{V})(\overline p)$ is an irreducible $L$-module.

{\bf Proof} Consider the map for a fixed $p \in 
\ZZ^n$ such that  ${\overline p} \neq 0$,
$$\pi(\overline p) : L( \stackrel{\sim}{V})(\overline p) \rightarrow
               L(\stackrel{\sim}{V})(\overline 0)$$ 
$$\pi(\overline p) ({\overline{v}_k} \otimes t^{k+r+p})= {\overline{v}_k} 
\otimes t^{k+r}.$$
It is easy to see that $\pi(\overline p)$ is a vector space isomorphism and 
not a  $L$-module map. For example 
$$(u, k+r+p +\alpha) {\overline{v}_k} \otimes t^{k+r+p} = 
d_i({\overline{v}_k} \otimes t^{k+r+p}) \neq d_i({\overline{v}_k} \otimes 
t^{k+r}) = (u, k+r +\alpha) {\overline{v}_k} \otimes t^{k+r}.$$
Now suppose $W$ is a non-zero proper submodule of $L( \stackrel{\sim}{V})
(\overline p)$. Clearly $\pi(\overline p)(W)$ is a non-zero proper subspace
of 
$L(\stackrel{\sim}{V})(\overline 0)$. To prove that $L( \stackrel{\sim}{V})
(\overline p)$ is irreducible, it is sufficient to prove that 
$\pi(\overline p)(W)$ is  an $L$-module.
Since $\pi(\overline p)$ commutes with 
$ L(\stackrel{\circ}{\gg}, {\bf \sigma}) \oplus  A(m) \oplus 
\displaystyle{\sum_{r \in \Gamma}} \CC t^r d_0$, it follows that 
$\pi(\overline p)(W)$  is a module for the above 
space. Since $W$ is a weight module, it is easy to check that 
$\pi(\overline p)(W)$ is also a weight module.   Let $w = {\overline{v}_k} 
\otimes t^{k+r} \in \pi(\overline p)(W)$ be a weight vector, then
$$D(u, s) ({\overline{v}_k} \otimes t^{k+r+p}) = 
\overline{I(u,s) v_k} \otimes t^{k+r+s+p} + (u, k+r+p +\alpha) 
{\overline{v}_k} \otimes t^{k+r+s+p} \in W.$$ Also
$$t^s( {\overline{v}_k} \otimes t^{k+r+p}) = {\overline{v}_k} 
\otimes t^{k+r+s+p} \in W.$$ This proves
$\overline{I(u,s) v_k} \otimes t^{k+r+s+p} + (u, k+r +\alpha) 
{\overline{v}_k} \otimes t^{k+r+s+p} \in W$.
This means $\pi^{-1}(\overline p)( D(u, s) ({\overline{v}_k} \otimes t^{k+r}))
 \in W$. So $D(u, s) ({\overline{v}_k} \otimes t^{k+r}) \in \pi(\overline p)
(W)$. Thus
$\pi(\overline p)(W)$  is an $L$-module. This proves that 
$L(\stackrel{\sim}{V})(\overline p)$ is irreducible $L$-module.

It is possible that some of the modules in $L(\stackrel{\sim}{V})$ are 
isomorphic. We need to develop the notion of graded automorphisms of 
$\stackrel{\sim}{V}$.
\item [{(8.11)}] {\bf Definition} An $\stackrel{\sim}{L}$-module 
automorphism 
$\theta$ of $\stackrel{\sim}{V}$ is called $\overline p$-graded if 
$\theta(\stackrel{\sim}{V}_{\overline k}) = \stackrel{\sim}{V}_{\overline k - 
\overline p}$ for all $\overline k \in \Lambda$.
\item [{(8.11.1)}] In this case dim $\stackrel{\sim}{V}_{\overline k}$ = dim 
$\stackrel{\sim}{V}_{\overline k - \overline p}$.
\item [{(8.11.2)}]Suppose $\theta$  is a $\overline p$-graded automorphism 
of $\stackrel{\sim}{V}$. Choose minimal integer $N_p > 0$ such that $N_p.
{\overline p} = 0$ in $\Lambda$.Then clearly $\theta^{N_p} 
\stackrel{\sim}{V}_{\overline k} = \stackrel{\sim}{V}_{\overline k}$.  Thus 
there exists a vector $v$ in $\stackrel{\sim}{V}_{\overline k}$ such that 
$\theta^{N_p} v = \lambda v$ for some  non-zero scalar $\lambda$.
Consider the space
$$W = \{ v \in \stackrel{\sim}{V} | \theta^{N_p} v = \lambda v \},$$
which can be seen as a graded submodule of $\stackrel{\sim}{V}$. Since 
$\stackrel{\sim}{V}$ is graded irreducible by Lemma 8.9, we see that $W = 
\stackrel{\sim}{V}$. Thus $\theta^{N_p} = \lambda$ on $\stackrel{\sim}{V}$.
Now by suitably multiplying  $\theta$  by a scalar we can assume 
$\theta^{N_p} = 1$.

So here after we will work only with graded automorphisms of finite order.
Recell we have fixed  $\alpha \in \CC^n$ from 8.4. For any $L$-module $V$, 
we define
$$V_k = \{ v \in V | D(u,0) v = (u, k + \alpha) v \},$$ for $k \in \ZZ^n$.
Thus we have 
$$L(\stackrel{\sim}{V})_k =  \displaystyle{\bigoplus _{{\overline p} \in \Lambda}
L( \stackrel{\sim}{V})(\overline p)_k}.$$
\item [{(8.12)}] {\bf Lemma } dim $ { \stackrel{\sim}{V}}_{\overline k} $ = 
dim $L( \stackrel{\sim}{V})(\overline 0)_k$ =   dim $L( \stackrel{\sim}{V})
(\overline p)_{k+p}$. \\
{\bf Proof} Proof follows from the definitions of $L( \stackrel{\sim}{V})
(\overline p)$.

\item [{(8.13)}] Suppose $\theta$  is a $\overline p$-graded automorphism of 
$\stackrel{\sim}{V}$ and $N_p$ be the order. Then $\theta$ defines an $L$
-module isomorphism
 $$\theta^{'} : L( \stackrel{\sim}{V})(\overline 0) \rightarrow  
L(\stackrel{\sim}{V})(\overline p)$$         
$$\theta^{'}({\overline{v}_k} \otimes t^{k+r}) = \theta({\overline{v}_k} ) 
\otimes t^{k+r+p}.$$

It is easy to check that $\theta^{'}$ is an  isomorphism of $L$-modules.
Suppose there exists an $L$-module isomorphism  $\theta^{'}$ from 
$L(\stackrel{\sim}{V})(\overline 0)$ to $L( \stackrel{\sim}{V})
(\overline p)$, then $k$-weight vectors go to $k$-weight vectors under this 
isomorphism. Thus 
$$\theta^{'}({\overline{v}_k} \otimes t^{k}) = {\overline{w}_{k-p}} \otimes 
t^{k}.$$
We now define $\theta({\overline{v}_k}) = {\overline{w}_{k-p}}$. This can be 
checked to  be $\overline p$-graded automorphism of 
$\stackrel{\sim}{V}$ and can be assumed to be of finite order.

\item [{(8.14)}] {\bf Proposition} There is a one-one correspondence between 
$\overline p$-graded $\stackrel{\sim}{L}$-module automorphism of 
$\stackrel{\sim}{V}$ and isomorphism between  $L( \stackrel{\sim}{V})
(\overline 0)$ and  $L( \stackrel{\sim}{V})(\overline p)$. \\
{\bf Proof} Proof follows from above discussion.

\item [{(8.15)}] {\bf Lemma} Suppose such an automorphism exists, then \\
 (1) dim $ { \stackrel{\sim}{V}}_{\overline k} $ = dim 
$ {\stackrel{\sim}{V}}_{\overline k - \overline p }$, \\
(2) dim $L( \stackrel{\sim}{V})(\overline 0)_k$ = dim 
$L(\stackrel{\sim}{V})(\overline 0)_{k-p}$. 

\item [{(8.16)}] {\bf Corollary} (1) dim 
${\stackrel{\sim}{V}}_{\overline k} $ 
 = dim $ { \stackrel{\sim}{V}}_{\overline {k + jp}}, j \in \ZZ $. \\
(2) dim $L({ \stackrel{\sim}{V}})(j \overline p )_{\overline {k + ip}}$ = dim 
$ {\stackrel{\sim}{V}}_{\overline k} $ for any $i,j \in \ZZ $. \\
{\bf Proof} Repeated application of $\theta$.  Note that $\theta^j$
is also an  automorphism.

\item [{(8.17)}] {\bf Theorem} $\stackrel{\sim}{V}$ is an $\stackrel{\sim}
{L}$ irreducible module if and only if $L( \stackrel{\sim}{V})(\overline p),
{ \overline p} \in \Lambda $ are mutually non-isomorphic as ${L}$-modules. \\
{\bf Proof } We can suppose $L( \stackrel{\sim}{V})(\overline 0) \cong 
L( \stackrel{\sim}{V})(\overline p)$ for $0 \neq {\overline p} \in \Lambda$.
 Then there exists a $\overline p$-graded automorphism of 
$\stackrel{\sim}{V}$. Let $N_p$ be the order. For ${\overline{v}_k} \in 
{\stackrel{\sim}{V}}_{\overline k} $, define 
$${\overline{v}_k}(0) = {\overline{v}_k} + \theta({\overline{v}_k}) + \cdots
 + \theta^{N_p -1}({\overline{v}_k}).$$
It is easy to check that  
$\theta({\overline{v}_k}(0)) = {\overline{v}_k}(0).$
Let \\
$\stackrel{\sim}{M}_0 = \{ {\overline{v}_k}(0), {\overline{v}_k} \in 
{\stackrel{\sim}{V}}_{\overline k}, { \overline k} \in \Lambda \} $, where $ 
\{  \} $ means the space generated by the vectors inside 
$\stackrel{\sim}{V}$. Clearly 
$\stackrel{\sim}{M}_0$ is a non-zero proper submodule of 
$\stackrel{\sim}{V}$.
This proves one side of the theorem.
Now suppose $L( \stackrel{\sim}{V})(\overline p), { \overline p}  \in 
\Lambda$ be mutually non-isomorphic modules. Suppose $W$ is a  
$\stackrel{\sim}{L}$ submodule of $\stackrel{\sim}{V}$. Then clearly
$$L(W) \subset L( \stackrel{\sim}{V}) = 
\displaystyle{\bigoplus _{{\overline p} \in \Lambda }L( \stackrel{\sim}{V})
(\overline p)}.$$
Then by Proposition 8.8,  we see that there exists $S \subset \Lambda $ 
such that
$$L(W)  = 
 \displaystyle{\bigoplus _{{\overline p} \in S }L( \stackrel{\sim}{V})
(\overline p)}.$$
This means $L( \stackrel{\sim}{V})(\overline p) \subset L(W)$ for some 
$\overline p$. This means ${\overline{v}_k} \otimes t^{k+p}
\in L(W)$, which means $\overline{v}_k \in W$. By Lemma 8.9, the module 
generated by $\overline{v}_k$ is $\stackrel{\sim}{V}$. Hence 
$W = \stackrel{\sim}{V}$. This proves the other part of the theorem.

\item [{(8.18)}] The aim is to find suitable irreducible 
$\stackrel{\sim}{L}$ submodule  of 
$\stackrel{\sim}{V}$.

Suppose $\stackrel{\sim}{V}$ is irreducible, then we are done. If not then, 
$L( \stackrel{\sim}{V})(\overline 0) \cong L( \stackrel{\sim}{V})(\overline 
p)$.
We will now find conditions for $\stackrel{\sim}{M}_0$  to be irreducible.
Let $\zeta$ be the $N_p$th primitive root of unity.
Define 
$${\overline{v}_k}(i) = {\overline{v}_k} +\zeta^{i} \theta({\overline{v}_k})
+ \cdots +\zeta^{(N_p - 1)i} \theta^{N_p -1}({\overline{v}_k}).$$
It is easy to check that 
$\theta({\overline v}_k(i)) = \zeta^{-i}  \ \ {\overline v}_k(i)$.
Let $$\stackrel{\sim}{M}_i = \{ v \in \stackrel{\sim}{V} | \theta(v) = 
\zeta^{-i} v \}, $$
which can be seen to be proper submodule of $\stackrel{\sim}{V}$.
Further
\item [{(8.19)}] $\stackrel{\sim}{V} = 
\displaystyle{\bigoplus_{i =0 }^{N_p -1} \stackrel{\sim}{M}_i},$
which is an eigen space decomposition.
\item [{(8.20)}] Define ${M}_i = \{{\overline{v}_k}(i) \otimes t^{k+r}, 
{\overline{v}_k} \in  { \stackrel{\sim}{V}}_{\overline k} , r \in \Gamma  \}.$
It can be verified to be an $L$-module and further
$$M_i \subset \displaystyle{\bigoplus_{j =1 }^{N_p -1} 
L( \stackrel{\sim}{V})(j p)}.$$

{\bf Claim} ${\overline{v}_k}(i) \otimes t^{k-jp+r} \in M_i$. \\
It is easy to check that
$$
{\overline{v}_k}(i) =  \zeta^i \theta({\overline{v}_k}(i)) \\
= \zeta^{ij} \theta^j({\overline{v}_k}(i)) \\
= {\overline{w}_{k-jp}(i)} .
$$
By definition ${\overline{w}_{k-jp}(i)} \otimes t^{k-jp+r} \in M_i$,
thus ${\overline{w}_{k-jp}(i)} \otimes t^{k-jp+r}  = {\overline{v}_k}(i) 
\otimes t^{k-j+r}  \in M_i$.
Hence the claim follows.
\item [{(8.20.1)}]In the Definition 8.20, one can replace $\Gamma$ by 
 $\Gamma_p =  \Gamma + \ZZ_p$. Note that $\ZZ_p$ is a finite subgroup of  
$\Lambda$.

\item [{(8.21)}] {\bf Lemma} 1)$M_i \cong   L( \stackrel{\sim}{V})( 0)$as 
$L$-module and in particular      $M_i$ is an irreducible $L$-module. \\
2) $\displaystyle{\bigoplus_{i =0 }^{N_p -1} M_i} = 
\displaystyle{\bigoplus_{j =0 }^{N_p -1} L( \stackrel{\sim}{V})(  j p)}$. \\
{\bf Proof} We know  that $L( \stackrel{\sim}{V})(\overline 0) \cong 
L(\stackrel{\sim}{V})(\overline { jp})$ and $M_i$ is contained in 
$\displaystyle{\bigoplus_{j =0 }^{N_p -1} 
L( \stackrel{\sim}{V})(j p)}$. Thus $M_i$ is isomorphic to the sum of 
finitely many copies of $L(\stackrel{\sim}{V})(\overline 0) $. We will now 
compare the dimensions of the weight spaces. \\
{\bf Claim} dim $ { \stackrel{\sim}{V}}_{\overline k}$  = dim $(M_i)_k$. \\
Consider the map ${\overline{v}_k} \rightarrow {\overline{v}_k}(i)$, which 
is clearly injective.  The surjectivity  is also obvious and hence  the 
claim.

>From Lemma 8.12, we know that dim $
 { \stackrel{\sim}{V}}_{\overline k}$ = dim $L( \stackrel{\sim}{V})(\overline
 0)_k$. This proves that the dimensions of the weight spaces of $M_i$ and 
$L( \stackrel{\sim}{V})(\overline 0)$ are same. Thus $M_i \cong 
L(\stackrel{\sim}{V})(\overline 0)$. In particular $M_i$ is an 
ireducible $L$-module. It is easy to see that the sum at LHS is direct. 
Equality holds for dimensions reason. Hence both parts of  the Lemma follows.
\item [{(8.22)}] Let $\ZZ_p$ be the cyclic group generated by $\overline p$ 
inside $\Lambda$. Let $\Lambda_p = \Lambda / \ZZ_p$.
Consider $$ \ZZ^n \rightarrow  \ZZ^n/ \Gamma = \Lambda \rightarrow \Lambda_p.
$$
Let $\Gamma_p$ be the kernal of the above map.
\item [{(8.22.1)}]Note that each 
$\stackrel{\sim}{M_i}$ is no more graded by $\Lambda$ but by $\Lambda_p$. 
Similar to the proof of Lemma 8.9, one can prove that each $\stackrel{\sim}
{M_i}$ is $\Lambda_p$-graded irreducible.

 Let $\overline 0 = q_0, q_1, ...,q_{n_p - 1}$ be a set of coset 
representative  for $\Lambda_p$. Note that  
$n_p N_p = |\Lambda |$. Since 
$\theta$  is a $\overline p$-graded automorphism, we see that 
$$L(\stackrel{\sim}{V})(q_i) \cong L(\stackrel{\sim}{V})(q_i + \overline p)
 \cdots \cong L(\stackrel{\sim}{V})(q_i + (N_p -1) \overline p).$$
Put $W(q_i) = \displaystyle{\bigoplus_{j =1 }^{N_p -1} 
L( \stackrel{\sim}{V})(q_i + j\overline   p)}$.  \\
Define for $0 \leq i < N_p,  \  \  0 \leq l < n_p$,
$$M_i^l = \{{\overline{v}_k}(i) \otimes t^{k+r +{q_l}}, {\overline{v}_k} \in 
{ \stackrel{\sim}{V}}_{\overline k} , r \in \Gamma  \}.$$
Note that $ M_i^{0} = M_i $.  \\
This can be verified to be an $L$-module and $M_i^{l} \subset W(q_l)$. \\
{\bf Claim} $M_i^{l}$ is irreducible $L$-module and isomorphic to 
$L(\stackrel{\sim}{V})(q_l)$. \\
The proof is similar to the case $l=0$. From the definition of  $W(q_l)$, 
we know that $M_i^{l}$ is isomorphic to sum of $L(\stackrel{\sim}{V})(q_l)$. 
As seen earliar for the case $l=0$, dim $(M_i^{l})_{k+q_l} \cong$ dim 
$ {\stackrel{\sim}{V}}_{\overline k} $, by considering the map
$v_k \rightarrow {{v}_k}(i) \otimes t^{k+q_l}$. 
But by Lemma 8.12, we know that 
dim 
$ {\stackrel{\sim}{V}}_{\overline k}$  = dim 
$L(\stackrel{\sim}{V})(q_l)_{k+q_l}$.
Thus
\item [{(8.23)}] $M_i^{l} \cong L( \stackrel{\sim}{V})(q_l),$
 and in particular $M_i^{l}$ is irreducible.

\item [{(8.24)}] {\bf Lemma} $L(\stackrel{\sim}{M_i}) = 
 \displaystyle{\bigoplus_{l=0 }^{n_p -1} M_i^{l}}$. \\
{\bf Proof }Proof follows from the Definition in (8.4). 
\item [{(8.25)}]We notice that R. H. S. in (8.23) is independent of $i$.
Thus $$L(\stackrel{\sim}{M_i}) \cong L(\stackrel{\sim}{M_j}).$$
Now we will record a simple fact.

\item [{(8.26)}]Let $W_1$ and $W_2$  be $\stackrel{\sim}{L}$-submodules of 
$\stackrel{\sim}{V}$. Then $W_1 \cong W_2$  as  $\stackrel{\sim}{L}$-modules
if and only if ${L}(W_1) \cong {L}(W_2)$ as $L$-modules. 
The above follows from the definitions. We now have the following.

\item [{(8.27)}] {\bf Proposition} $\stackrel{\sim}{M_i} \cong 
\stackrel{\sim}{M_j}$. \\
{\bf Proof }Proof follows from (8.26). \\
Now we have $${L}(\stackrel{\sim}{M_i}) = 
\displaystyle{\bigoplus_{l=0 }^{n_p -1} M_i^{l}},$$
where each $M_i^{l}$ is irreducible  $L$-module. This situation is similar 
to Theorem 8.17. We can now prove the following, whose proof is similar to 
Theorem 8.17.
\item [{(8.28)}] {\bf Theorem}
We fix $i$. $\stackrel{\sim}{M_i}$ is irreducible as $\stackrel{\sim}{L}$
-module  if and only if  $M_i^{l}, 0 \leq l < {n_p} $ are mutually 
non-isomorphic as $L$-modules.

{\bf Proof of the Theorem (8.7)}
Recall that $\stackrel{\sim}{V} = 
\displaystyle{\bigoplus_{i =0 }^{N_p -1} \stackrel{\sim}{M_i}}$ and all 
components are isomorphic as $\stackrel{\sim}{L}$-modules, which follows 
from $(8.19)$ and 
$(8.27)$. Further dim $\stackrel{\sim}{M_i} < $  dim $\stackrel{\sim}{V}$. 
Suppose 
$\stackrel{\sim}{M_i}$ is reducible as $\stackrel{\sim}{L}$-module, then we 
can repeat the above process. This process has to stop for dimensions 
reasons. Hence the theorem follows.

\item [{(8.29)}] To avoid more complex notation we assume that each 
$\stackrel{\sim}{M_i}$ is irreducible  $\stackrel{\sim}{L}$-module.

\end{enumerate}

\section{Final Theorem}

In this section we will describe the $\stackrel{\sim}{L}$-module structure
 of $\stackrel{\sim}{M_i}$. We are assuming  that each
$\stackrel{\sim}{M_i}$ is irreducible     $\stackrel{\sim}{L}$-module.
We will actually prove that $\stackrel{\sim}{M_i}$ is an irreducible module 
for the direct sum $g l_n \oplus \stackrel{\circ}{\gg}$. Recall that 
$gl_n$ is a quotient of $I$ from 7.11.3 and 
$\stackrel{\circ}{\gg}$ is a quotient of $L(\stackrel{\circ}{\gg}, \sigma)$ 
by the map $X(k) \rightarrow X \in\stackrel{\circ}{\gg}_{\overline{k}}$.  Here 
$X(\overline k)$ is identified by $X$.We will start with a simple Lemma. 

\quad Suppose $S$ is a vector space such that $S =
\displaystyle{\bigoplus_{\overline{k} \in \Lambda}} S_{\overline{k}}$. 

An operator  $T:S \rightarrow S$ is called degree $\overline{k}$ 
operator 
if $T( S_{\overline{l}}) \subseteq  S_{\overline{l} + \overline{k}} $. The 
following lemma is trivial to see.

\begin{enumerate}

\item [{(9.1)}] {\bf Lemma} Suppose $T$ is a degree $\overline{k}$ operator 
such that $\overline{k} \neq 0$. 
Suppose $T$ acts as a scalar $\lambda$. Then $\lambda=0$. 

\quad The following is well known. See Proposition 19.1(b) of $[H]$.

\item [{(9.2)}] {\bf Lemma} Let $\gg^\prime$ be a Lie algebra which need not 
be finite dimensional. Let 
$(V_1, \rho)$ be an irreducible finite dimensional module for $\gg^\prime$. 
We have a map 
$\rho : \gg^\prime \rightarrow End \ V_1$. Then $\rho(\gg^\prime)$ is a 
reductive Lie-algebra with at most one dimensional center.

\quad Let $\gg_1$ and $\gg_2$ be infinite dimensional Lie algebra such that 
$\gg_1$ acts on $\gg_2$. Let 
$\gg^\prime = \gg_1 \ltimes \gg_2$ be the natural semi-direct Lie algebra. 
Let  $J$ be an abelian ideal of $\gg_2$ which will not be assumed to be a
 ideal of $\gg^\prime$.

\item [{(9.3)}] {\bf Proposition} Suppose $(V^\prime, \rho)$ is an 
irreducible finite dimensional module for 
$\gg^\prime$. We have $\rho : \gg^\prime \rightarrow End (V^\prime)$. Then 
$\rho(J)$ is central ideal in $\rho(\gg^\prime).$\\
{\bf Proof} From above Lemma (9.2), it follows that $\rho(\gg^\prime)$ is a 
reductive Lie algebra. Since $\gg_2$ 
is an ideal in $\gg^\prime$ we have $\rho(\gg_2)$ is an ideal in 
$\rho(\gg^\prime)$. Thus $\rho(\gg_2)$ 
is reductive Lie algebra. Now we know that $J$ is abelian ideal in $\gg_2$ 
and hence $\rho(J)$ is contained in the center of $\rho(\gg_2)$. This proves
$\rho(J)$ is actually contained in the center of $\rho(\gg^\prime)$. In 
particular $\rho(J)$ is an ideal in $\rho(\gg^\prime)$.
\item [{(9.4)}] {\bf Theorem} $\stackrel{\sim}{M_i}$ is actually an 
irreducible finite dimensional module for the 
direct sum $g l_n \oplus \stackrel{\circ}{\gg}$.\\
We need the following Lemma. 
\item [{(9.4.1)}] {\bf Lemma}Let $k \in \ZZ^n, r_1, r_2, \cdots, r_d \in 
\Gamma$ and $d \geq 1$. Suppose $X(k, r_1, \cdots, r_d)$ acts as a scalar 
$\lambda$ on $\stackrel{\sim}{V}$. Then the scalar $\lambda$ is zero. \\
{\bf Proof} First we note that if  $X \in \stackrel{\circ}{\gg} \cap 
\gg(\overline{\circ}, \overline{\circ})$ then 
$X \in \hh(0)$.  This follows from  the theory of finite dimensional simple 
Lie algebras. Suppose $\overline{k}=0$, then $X \in \hh(0)$.  From 6.3.5, 
it follows that $X(k,r) = X(k) - X(k+r)$ is zero on $\stackrel{\sim}{V}$. 
Since $X(k,r_1, r_2, \cdots, r_d)$ is spanned by   $X(k,r), r \in \Gamma$,
it follows that $\lambda$ is zero. Now suppose $\overline{k} \neq 0$, then 
from Lemma 9.1, it follows that $\lambda$ is zero. This completes the proof 
of the Lemma.

{\bf Proof of Theorem 9.4} 
Let $\rho$  denote the  $\stackrel{\sim}{L}$-module action on 
$\stackrel{\sim}{M_i}$ and note that  $\rho$  is independent of $i$ as all 
$\stackrel{\sim}{M_i}$ are isomorphic. Consider $(\ker \rho) \cap I$ 
which is a co-finite ideal of $I$. Thus from 7.11.2, it follows that 
$(\ker \rho) \cap I$ contains $I_d$ for large $d$. Consider 
$$[I(u,k, r_1, r_2, \cdots, r_d), X(l)] = X(k+l, r_1, r_2, \cdots, r_d) \in 
\ker \rho. \ \ (*)$$ 
This  proves $F_d$  acts trivially on $\stackrel{\sim}{M_i}$. But we know 
that $[F_{d-1}, F_{d-1}] \subseteq F_d$ by 7.8.2. Thus   $\rho(F_{d-1})$ is an
abelian ideal in  $\rho(L(\stackrel{\circ}{\gg}, \sigma))$. By Proposition 
9.3 it follows that $\rho(F_{d-1})$ is a central ideal in 
$\rho(\stackrel{\sim}{L})$. It is well known that center acts as  scalars 
on a finite dimensional irreducible module. Thus  
$X(k,r_1,r_2, \cdots, r_d)$ acts  as scalar on each
$\stackrel{\sim}{M_i}$ and the scalar is independent of $i$. 
Thus  $X(k,r_1,r_2, \cdots, r_d)$ acts  as  a single scalar on 
$\stackrel{\sim}{V}$. Now by Lemma 9.4.1, scalar is zero. 
Thus  $F_{d-1}$ acts trivially on each
$\stackrel{\sim}{M_i}$. It now follows that $\rho(I_{d-1})$ is an ideal in 
$\rho(\stackrel{\sim}{L})$ by  $(*)$. By Lemma 4.2 of [E2], $I_{d-1}$ acts 
trivially. By repeating this argument finitely many times we see that  
$I_2$ acts trivially and $F_2$ acts trivially. Now from above argument we 
see that $F_1$ acts trivially. Thus  $I_2 \oplus F_1$ acts trivially on 
$\stackrel{\sim}{M_i}$. As noted in the beginning of the section we see that
each $\stackrel{\sim}{M_i}$ is a module for $g l_n \oplus 
\stackrel{\circ}{\gg}$ and hence the theorem is proved.

\item [{(9.5)}] By Lemma 2.7 of $[HL]$, there exists an irreducible module 
$\stackrel{\sim}{V}_1$ of $gl_n$ and an irreducible module 
$\stackrel{\sim}{V}_2$ for $\stackrel{\circ}{\gg}$ such that 
$\stackrel{\sim}{V}_1 \otimes \stackrel{\sim}{V}_2 \cong \stackrel{\sim}{V}$
as $g l_n \oplus \stackrel{\circ}{\gg}$-module.
Regarding gradation, note that all vectors of $I$ are grade zero and hence 
we can assume $\stackrel{\sim}{V}_1$  
is zero graded and $\stackrel{\sim}{V}_2$  is $\Lambda_p$-graded.

\item [{(9.6)}] We will now describe $L$-module $T$ in terms of $\stackrel
{\sim}{V}_1$ and $\stackrel{\sim}{V}_2$.\\
\quad Let $\{E_{ij}\}_{1 \leq i, j \leq n}$ be the standard basis of $gl_n$. 
From $[E2]$ it follows that $I(u,r) = \sum u_i r_j m_j E_{ji} \in gl_n$ is 
linear in both variables mod $I_2$.\\
Let $u=\sum u_i e_i$ and $r= \sum m_j r_j e_j$,\\
then $D(u,r)= D(u,0) + \displaystyle{\sum_{i,j}} u_i r_j m_j E_{ji}$  mod 
$I_2$.\\
Let $\stackrel{\sim}{V}_2= \displaystyle{\bigoplus_{\overline{k} \in \Lambda_p}}
\stackrel{\sim}{V}_{2, \overline{k}}$ 
be the $\Lambda_p$-gradation.\\
We will now define $L$-module on $\stackrel{\sim}{V}_1 \otimes \stackrel
{\sim}{V}_2 \otimes A_n$.
$$
\begin{array}{lll}
D(u,r) (v_1 \otimes v_2 \otimes t^k) &= (u, k+\alpha) v_1 \otimes v_2 
\otimes t^{k+r}\\
&+\sum(u_i r_jm_j E_{ji} v_1) \otimes v_2 \otimes t^{k+r},\\
X(l) (v_1 \otimes v_2 \otimes t^k) &= v_1 \otimes X  v_2 \otimes t^{k+l},\\
t^r \cdot v_1 \otimes v_2 \otimes t^k &= v_1 \otimes  v_2 \otimes t^{k+r},\\
t^r d_0 \cdot v_1 \otimes v_2 \otimes t^k &= \lambda(d_0) v_1 \otimes v_2 
\otimes t^{k+r}.\\
\end{array}
$$

From the above discussion we see that  this module is isomorphic to 
$L(\stackrel{\sim}{M_i})$ as $L$-module.\\
Let $\displaystyle{\bigoplus_{\underline{k} \in \Lambda_p}}
\stackrel{\sim}{V}_{2, \underline{k}}= \stackrel{\sim}{V}_2$ 
be the $\Lambda_p$-gradation.\\ 
Then consider the submodule of $L(\stackrel{\sim}{M_i})$\\
$\displaystyle{\bigoplus_{k \in \ZZ^n}}\stackrel{\sim}{V}_1 \otimes 
\stackrel{\sim}{V}_{2, {\underline k}} \otimes t^k$ which is easy 
to see that it is isomorphic to $M_i$  for any $i$ as $L$-module and  hence 
irreducible $L$-module. By defining  $t^r K_p$ to be zero for $1 \leq p 
\leq n$ we see that the above module is actually a $\tau^0$-module.
\item [{(9.6.1)}] We have seen that $T \cong 
L(\stackrel{\sim}{V})(0)$ as $L$-modules from (8.6). But $M_i \cong 
L(\stackrel{\sim}{V})(0)$ by 8.21. Thus it follows that $M_i \cong T$ as 
$L$-modules for each $i$. Note that it is difficult to give direct module
map between $M_i$ and  $T$ as there is a twisting taking place. Thus we 
described $T$ explicitly as $L$-module.
Let $M= Ind_{\tau_0+\tau^+}^\tau T$\\
Then there exists a unique maximal submodule $M^{rad}$ intersecting $T$ 
trivially. Thus $M/M^{rad}$ is irreducible and 
isomorphic to the original module $V$.

\item [{(9.7)}] {\bf Theorem} Let $V$ be an irreducible integrable module 
for $\tau$ with $K_0$ acts as 
$C_0 > 0$ and $K_i$ acts trivially. Let $T$ and $M$ as above. Then $V \cong 
M/M^{rad}$ as $\tau$-modules.
\end{enumerate}

{\bf Acknowledgements:}\

\quad The first author is very grateful to Professor Kaiming Zhao for his 
kind invitation to Chinese Academy of Sciences, China where some of the 
work was done.


\begin{thebibliography}{BMS}

\bibitem[ABFP]{ABFP} Allison, B.,  Berman, S., Faulkner, J and Pianzola, 
\emph{A. Multiloop realization of Extended Affine Lie algebras and the Lie 
Tori,} Trans. Amer. Math. Soc., {\bf 361}(2009), 4807-4842.

\bibitem[BB]{BB} Berman, S and Billig, 
\emph{Y. Irreducible representations for toroidal Lie algebras,}
Journal of Algebra, {\bf 221}(1999), 188-231.

\bibitem[B]{B} Billig, Y. \emph{Jet Modules, Canad.}
Journal Math. {\bf 59(4)}(2007), 712-729.

\bibitem[B2]{B2} Billig, Y and Michael Lau, \emph{Thin covering of modules,}
Journal of Algebra, {\bf 316}(2007), 147-173.

\bibitem[E1]{E1} Eswara Rao, 
\emph{S. Irreducible representations of the Lie algebra of the 
diffeomorphisms of a d-dimensional torus,} Journal of Algebra, 
{\bf 182}(1996), 401-421. \bibitem[E2]{E2} Eswara Rao, 
\emph{S. Partial classification of modules for Lie algebra of 
diffeomorphisms of d-dimensional torus, }
Journal of Math. Physics, {\bf 45}(8)(2004), 3322-3333.
\bibitem[E3]{E3} Eswara Rao, \emph{S. Complete reducibility for the 
twisted affine Lie algebras,}
Communications in Algebra, {\bf 40}(2)(2012), 379-385.

\bibitem[E4]{E4} Eswara Rao, S. and Sachin Sharma, \emph{ Integrable 
modules for Lie Torus,}To be published in Journal of Pure and Applied 
Algebra.

\bibitem[EJ]{EJ} Eswara Rao, S. and Cuipo Jiang, 
\emph{Classifications of irreducible integrable representations for the 
full toroidal Lie algebra,}
Journal of Pure and Applied Algebra, {\bf 200}(2005), 71-85.

\bibitem[EMY]{EMY} Eswara Rao, S., Moody, R.V. and Yokonuma,
\emph{T. Toroidal Lie algebras and vertex representations,}
Geom, Dedicata, {\bf 35}(1990), 283-307.

\bibitem[FJ]{FJ} Fu Jiayuan and Cuipo Jiang,
\emph{Integrable representations for the twisted full toroidal Lie algebra,}
Journal of Algebra {\bf 307}(2007), 769-794.

\bibitem[H]{H} Humphreys, \emph{J.E. Introduction to Lie algebras and 
Representation theory,} Springer, Berlin, Heidelberg and New York, 
{\bf 1972}.

\bibitem[HL]{HL} Haisheng Li,\emph{On certain categories of modules for 
Affine Lie algebra,} Mathematische zeitschrift, {\bf 248}(2004), 635-664.

\bibitem[FM]{FM} Jiang Cuipo and Meng Daoji, 
\emph{Integrable representations for gebaralized Virasoro-toroidal Lie 
algebra,} Journal of Algebra, {\bf 270}(2003), 307-334.

\bibitem[JS]{JS} Jie Sun, 
\emph{Universal central extensions of twisted forms of split simple Lie 
algebra over rings,} Journal of Algebra, {\bf 322}(2009), 1819-1829.

\bibitem[K]{K} KCa, V.G., \emph{Infinite dimensional Lie algebras,}
3rd ed. Cambridge University press, {\bf 1990}.

\bibitem[Ka]{Ka} Kassel, C., \emph{Kahler differentials and coverings 
extended over a commutative algebra,} Journal of Pure and Applied Algebra, 
{\bf 34}(1984), 265-275.

\bibitem[KN]{KN} Katsuyuki Naoi,
\emph{Multiloop Lie algebras and the construction of extended affine Lie 
algebra,} Journal of Algebra {\bf 323}(2010), 2103-2129.

\bibitem[L]{L} Lau Michael, \emph{Representation of Multiloop algebras,}
Pacific Journal of Mathematics, {\bf 245}(2010), 167-184.
\end{thebibliography}
\end{document}